
\magnification=1000
  
 \def\oui{oui} 
  
\ifx\arXiv\oui
\else
 \pdfpagewidth=210truemm
 \pdfpageheight=297truemm 
\fi
  
  %
  %

  %
  %

  \catcode`@=12 

 \def\defrefnote#1{\definexref{#1}{{\the\footnotenumber}}{refnotes}}

  %
  %


\ifx\couleurs\oui
\input graphicx
 \pdfpagewidth=210truemm
 \pdfpageheight=297truemm 
 \voffset=-5mm
\fi

\input eplain.tex
\expandafter\def\expandafter\newdimen\expandafter{\newdimen}

\ifx\couleurs\oui
\beginpackages
\usepackage{color}
\endpackages 
 \pdfpagewidth=210truemm
 \pdfpageheight=297truemm 
\long\def\rge#1{{\color{red}#1}}

\definecolor{bleu-iecn}{cmyk}{.98,.13,.1,.55}

\else
\long\def\rge#1{#1}

\fi

\makeatletter
\def\numberedfootnote{%
ÊÊ\global\advance\footnotenumber by 1
ÊÊ\@eplainfootnote{{\number\footnotenumber}}%
}%
\def\makecolumns#1/#2 {\par \begingroup
ÊÊ \@columndepth = #1
ÊÊ \advance\@columndepth by -1
ÊÊ \divide \@columndepth by #2
ÊÊ \advance\@columndepth by 1
ÊÊ \@linestogoincolumn = \@columndepth
ÊÊ \@linestogo = #1
ÊÊ \currentcolumn = 1
ÊÊ \def\@endcolumnactions{%
ÊÊÊÊÊÊ\ifnum \@linestogo<2
ÊÊÊÊÊÊÊÊ \the\crtok \egroup \endgroup \par 
ÊÊÊÊÊÊ\else
ÊÊÊÊÊÊÊÊ \global\advance\@linestogo by -1
ÊÊÊÊÊÊÊÊ \ifnum\@linestogoincolumn<2
ÊÊÊÊÊÊÊÊÊÊÊÊ\global\advance\currentcolumn by 1
ÊÊÊÊÊÊÊÊÊÊÊÊ\global\@linestogoincolumn = \@columndepth
ÊÊÊÊÊÊÊÊÊÊÊÊ\the\crtok
ÊÊÊÊÊÊÊÊ \else
ÊÊÊÊÊÊÊÊÊÊÊÊ&\global\advance\@linestogoincolumn by -1
ÊÊÊÊÊÊÊÊ \fi
ÊÊÊÊÊÊ\fi
ÊÊ }%
ÊÊ \makeactive\^^M
ÊÊ \letreturn \@endcolumnactions
ÊÊ \@columnwidth = \hsize
ÊÊÊÊ \advance\@columnwidth by -\parindent
ÊÊÊÊ \divide\@columnwidth by #2
ÊÊ \penalty\abovecolumnspenalty
ÊÊ \noindent 
ÊÊ \valign\bgroup
ÊÊÊÊ &\hbox to \@columnwidth{\strut \hsize = \@columnwidth ##\hfil}\cr
}%
\makeatother

\lefteqnumbers
   \def\testd{oui}
   \def\choixlat{\ifx\numadroite\testd\righteqnumbers
            \else  \lefteqnumbers\fi}
    \choixlat

\catcode`@=\letter
\def\@eplainfootnote#1{\let\@sf\empty 
  \ifhmode\edef\@sf{\spacefactor\the\spacefactor}\/\fi
  \global\advance\hlfootlabelnumber by 1
  \hlstart@impl{foot}{\hlfootlabel}%
  \hldest@impl{footback}{\hlfootbacklabel}%
  \hbox{$^{(#1)}$}%
  \hlend@impl{foot}%
  \@sf\vfootnote{#1.}%
}%
\catcode`@=\other

  \interfootnoteskip=0pt
  \let\note=\numberedfootnote
  \everyfootnote={\eightpoint\leftskip=5truemm\rightskip5truemm}
  
  \hsize150truemm\vsize 240truemm\hoffset=5truemm

  \pretolerance=500\tolerance=1000\brokenpenalty=5000
  \parindent3mm
  
  \countdef\temps=170
  \temps=\time
  \countdef\nminutes=171{\nminutes = \time}
  \countdef\nheure=172
  \def\heure{\begingroup                     
     \temps = \time \divide\temps by 60
     \nheure = \temps                        
     \nminutes = \time
     \multiply\temps by 60
     \advance\nminutes by -\temps            
     \ifnum\nminutes<10 \toks1 = {0}%
     \else\toks1 = {}%
     \fi
     \number\nheure h\the\toks1 \number\nminutes  
  \endgroup}%

  \newcount\chstart
  \chstart=\pageno
 \headline={\ifnum\pageno=\chstart {\hfill} \else {\hss \tenrm --\ \folio\ --\hss}\fi}
  \footline={\hfill}
  \normalbaselines
  \frenchspacing
    \def\dater{\vglue-10mm\rightline{(\the\day/\the\month/\the\year)}}
  \def\dateheure{\vglue-10mm\rightline{(\the\day/\the\month/\the\year,\ \heure)}}

  \newif\ifpagetitre \pagetitretrue
\newtoks\hautpagetitre \hautpagetitre={\hfill}
\newtoks\baspagetitre \baspagetitre={\hfill}
\newtoks\auteurcourant \auteurcourant={\hfill}
\newtoks\titrecourant \titrecourant={\hfill}
\newtoks\hautpagegauche
\newtoks\hautpagedroite
\newtoks\hautpagemilieu
\hautpagemilieu={\tenrm\hfil -- \folio\ -- \hfil}
\hautpagegauche={\ifx\midfolio\oui\the\hautpagemilieu\else\tenrm\folio\hfill\the\auteurcourant\hfill\fi}
\hautpagedroite={\ifx\midfolio\oui\the\hautpagemilieu\else\hfill\the\titrecourant\hfill\tenrm\folio\fi}
\newtoks\baspagegauche \baspagegauche={\hfil}
\newtoks\baspagedroite \baspagedroite={\hfil}
\headline={\ifpagetitre\the\hautpagetitre
\else\ifodd\pageno\the\hautpagedroite\else\the\hautpagegauche\fi\fi }
\footline={\ifpagetitre\the\baspagetitre
\else\ifodd\pageno\the\baspagedroite
\else\the\baspagegauche\fi\fi \global\pagetitrefalse}

\def\pageblanche{\vfill\eject\pagetitretrue
\null\vfill\eject
\pagetitretrue
}
\def\chgtpage{\ifodd\pageno \else
\pageblanche \fi \pagetitretrue\titreun=0\footnotenumber=0}

\def\chgtpageincrtitreun{\ifodd\pageno \else
\pageblanche \fi \pagetitretrue\footnotenumber=0}

\def\majnombres{\ifodd\pageno \else
\pageblanche \fi \pagetitretrue\hautpoly\titreun=0\footnotenumber=0}

\def\hautspages#1#2{\auteurcourant={\ninepcap#1}\titrecourant={\nineit#2}}

\ifnum\chstart=\pageno \pagetitretrue\fi
  

  \def\PAR{\par}
  
  \def\leftnote#1{\vadjust{\setbox1=\vtop{\hsize 20mm\parindent=0pt\eightpoint
  \baselineskip=9pt\rightskip=4mm plus 4mm\vskip-4.7mm#1}\hbox{\kern-2cm\smash{\box1}}}}

  
  \def\raggedcenter{\leftskip=20pt plus 10em  
       \rightskip=\leftskip 
        \parfillskip=0pt 
         \spaceskip=.3333em \xspaceskip=.5em 
          \pretolerance=9999 \tolerance=9999
           \hyphenpenalty=9999 \exhyphenpenalty=9999 }
           
  \def\titrecentre#1{{\parindent0mm\raggedcenter
       \spaceskip=.6em plus .2em minus .2em\xspaceskip=.6em plus .2em minus .2em
        \tit#1\par}}
        


  \def\oui{oui}
  
\def\fontetitreun{\ifx\paradouze\oui\douzepts\gpdouze\twelvebf\textfont1=\twelveib\else
\quatorzepts\gpquatorze\fourteenbf\fi}

\def\fontetitreunl{\douzepts\textfont1=\twelveib\scriptfont1=\tenib\fourteenti}
 
 \def\fontetitredeux{\textfont1=\eleveni\ifx\paradouze\oui\onzepts\scriptfont1=\ninei\elevenit\else
                        \douzepts\twelveit\fi}
 
   \def\fontetitredeuxb{\ifx\paradouze\oui\onzepts\eleventi\gponze\textfont1=\elevenib\scriptfont1=\nineib
                         \else\douzepts\twelveti\scriptfont1=\twelveib\scriptfont1=\tenib\gpdouze\fi}
                         
\def\fontetitredeuxl{\onzepts\textfont1=\elevenbf\scriptfont1=\ninebf\twelvebf}
  
\def\fontetitretrois{\textfont0=\elevenrm\scriptfont0=\eightrm\textfont1=\eleveni
                      \scriptfont1=\eighti\scriptscriptfont1=\sixi\elevenit}
                      
\def\fontetitrequatre{\textfont0=\elevenrm\scriptfont0=\eightrm\textfont1=\eleveni
                      \scriptfont1=\eighti\scriptscriptfont1=\sixi\elevenrm}
  
  \newcount\titreun\titreun=0
  \newcount\titredeux\titredeux=0
  \newcount\titretrois\titretrois=0
  \newcount\titrequatre\titrequatre=0
  \newcount\enonce\enonce=0
  
  \def\incr#1{\global\advance#1 by 1 {\the #1}}
  \def\avance#1{\global\advance#1 by 1}
  \def\init#1{\global#1=0}
  
  \long\def\Indentation#1#2{\setbox10=\hbox{\fontetitreun#1}
  	                    \ifdim\wd10 < 4mm
                         \setbox10=\hbox to 4mm{\box10\hfill}
                       \else\ifdim\wd10 < 6mm
                         \setbox10=\hbox to 6mm{\box10\hfill}
  	                    \else\ifdim\wd10 < 8mm
                         \setbox10=\hbox to 8mm{\box10\hfill}
                       \else\ifdim\wd10 < 12mm
                         \setbox10=\hbox to 12mm{\box10\hfill}
                       \fi\fi\fi\fi
                       \dimen10=\hsize
                       \advance \dimen10 by -\wd10
                       \noindent \box10 %
                       \ignorespaces
                       \hbox{\vtop{\hsize=\dimen10\raggedright\noindent\fontetitreun#2}}}

  \long\def\paraun#1{\removelastskip\par\medskip\goodbreak\vskip0pt plus.01\vsize\penalty-100
                \vskip0pt plus-.01\vsize
  	              \init{\titredeux}\ifnum\optionparag=1{\init\eqnumber\init\enonce}\else{}\fi
                  \goodbreak{\fontetitreun
  	                \Indentation{\incr{\titreun}.\ }{\fontetitreun #1\par}}\nobreak\medskip}

 %
 %
 \long\def\paraunc#1{\removelastskip\par\bigskip\goodbreak\vskip0pt plus.01\vsize\penalty-100
                \vskip0pt plus-.01\vsize
  	              \init{\titredeux}
                 \ifnum\optionparag=1{\init{\eqnumber}\init\enonce}\else{}\fi
                  \goodbreak
  	                {\parindent0mm\raggedcenter\fontetitreun\incr{\titreun}.\ 
                     \fontetitreun #1\par}\nobreak\medskip}
                     
\newtoks\titreunl
\titreunl={\ifnum\titreun=1{I}\fi%
\ifnum\titreun=2{II}\fi%
\ifnum\titreun=3{III}\fi%
\ifnum\titreun=4{IV}\fi%
\ifnum\titreun=5{V}\fi%
\ifnum\titreun=6{VI}\fi%
\ifnum\titreun=7{VII}\fi%
\ifnum\titreun=8{VIII}\fi%
\ifnum\titreun=9{IX}\fi%
\ifnum\titreun=10{X}\fi%
\ifnum\titreun=11{XI}\fi%
\ifnum\titreun=12{XII}\fi%
\ifnum\titreun=13{XIII}\fi%
}
\long\def\paraunl#1{\removelastskip\par\bigskip\bigskip\goodbreak\vskip0pt plus.01\vsize\penalty-100
                \vskip0pt plus-.01\vsize
  	              \init{\titredeux}\ifnum\optionparag=1{\init\eqnumber\init\enonce}\else{}\fi
                  \goodbreak{\fontetitreunl
  	                \Indentation{\global\advance\titreun by 1{\the\titreunl}.\ }{\fontetitreunl #1\par}}\nobreak\smallskip}

  
  \long\def\paradeux#1{\init{\titretrois}\vskip0pt plus.01\vsize\penalty-10
                \vskip0pt plus-.01\vsize\ifx \elie\oui\medskip\ifnum\titredeux>0\medskip\fi\fi
                 \Indentation{\fontetitredeux\the\titreun${\cdot}$\incr{\titredeux}.}
                              {\fontetitredeux\textfont1=\eleveni#1}\nobreak\par }
  
  \long\def\paradeuxb#1{\init{\titretrois}\vskip0pt plus.001\vsize\penalty-10
                \vskip0pt plus-.01\vsize{\ifx \elie\oui\medskip\ifnum\titredeux>0\medskip\fi\fi
                  \Indentation
  {\fontetitredeuxb\the\titreun${\cdot}$\incr{\titredeux}.}{ \fontetitredeuxb#1}}\nobreak
\smallskip}

\newtoks\titredeuxl
\titredeuxl={\ifnum\titredeux=1{A}\fi%
\ifnum\titredeux=2{B}\fi%
\ifnum\titredeux=3{C}\fi%
\ifnum\titredeux=4{D}\fi%
\ifnum\titredeux=5{E}\fi%
\ifnum\titredeux=6{F}\fi%
\ifnum\titredeux=7{G}\fi%
\ifnum\titredeux=8{H}\fi%
\ifnum\titredeux=9{I}\fi%
\ifnum\titredeux=10{J}\fi%
\ifnum\titredeux=11{K}\fi%
\ifnum\titredeux=12{L}\fi%
\ifnum\titredeux=13{M}\fi%
}
 \long\def\paradeuxl#1{\init{\titretrois}\vskip0pt plus.001\vsize\penalty-10
                \vskip0pt plus-.01
                \vsize \bigskip%
                  \Indentation
     {\fontetitredeuxl\global\advance\titredeux by 1
  \quad \the\titreunl${\cdot}$\the\titredeuxl.}{ \fontetitredeuxl#1}
  \removelastskip\nobreak\smallskip}
  

  \long\def\paratrois#1{\init{\titrequatre}\ifdim\lastskip<\smallskipamount
                \removelastskip\smallskip\fi
                 \vskip0pt plus.01\vsize\penalty-10
                  \vskip0pt
plus-.01\vsize{\ifx \elie\oui\ifnum\titretrois>0\medskip\fi\fi
\Indentation{\fontetitretrois\the\titreun${\cdot}$\the\titredeux${\cdot}$\incr{\titretrois}.\ }
  {\hskip0mm\baselineskip=14pt\fontetitretrois#1}\nobreak\smallskip}}
  
  
  \long\def\paratroisl#1{\init{\titrequatre}\ifdim\lastskip<\smallskipamount
                \removelastskip\fi
                 \vskip0pt plus.01\vsize\penalty-10
                  \vskip0pt
plus-.01\vsize\ifx \elie\oui\bigskip
\fi
\Indentation{\fontetitretrois\quad \quad \the\titreunl{${\cdot}$}\the\titredeuxl${\cdot}$\incr{\titretrois}.\ }
  {\hskip0mm\fontetitretrois#1}\nobreak\smallskip}


  \long\def\paraquatre#1{\ifdim\lastskip<\smallskipamount
                \removelastskip\smallskip\fi
                 \vskip0pt plus.01\vsize\penalty-10
                  \vskip0pt
                  plus-.01\vsize\par
 
\Indentation{\fontetitrequatre \the\titreun{${\cdot}$}\the\titredeux${\cdot}$\the\titretrois${\cdot}$\incr{\titrequatre}.\ }
{\hskip0mm\fontetitrequatre#1}\nobreak\smallskip}


\newtoks\titrequatrel
\titrequatrel={\ifnum\titrequatre=1{a}\fi%
\ifnum\titrequatre=2{b}\fi%
\ifnum\titrequatre=3{c}\fi%
\ifnum\titrequatre=4{d}\fi%
\ifnum\titrequatre=5{e}\fi%
\ifnum\titrequatre=6{f}\fi%
\ifnum\titrequatre=7{g}\fi%
\ifnum\titrequatre=8{h}\fi%
\ifnum\titrequatre=9{i}\fi%
\ifnum\titrequatre=10{j}\fi%
\ifnum\titrequatre=11{k}\fi%
\ifnum\titrequatre=12{l}\fi%
\ifnum\titrequatre=13{m}\fi%
}
\long\def\paraquatrel#1{\ifdim\lastskip<\smallskipamount
                \removelastskip\smallskip\fi
                 \vskip0pt plus.01\vsize\penalty-10
                  \vskip0pt
                  plus-.01\vsize{\bigskip
\Indentation{\global\advance\titrequatre by 1
\fontetitrequatre\quad \quad \quad \the\titreunl${\cdot}$\the\titredeuxl${\cdot}$\the\titretrois${\cdot}$\the\titrequatrel.\ }
{\hskip0mm\fontetitrequatre#1}\nobreak\smallskip}}

\ifx\optionkeys\oui
\def\drefun#1{\definexref{¤#1}{{\the\titreun}}{}} 
\def\drefdeux#1{\definexref{¤#1}{{\the\titreun}.{\the\titredeux}}{}}
\def\dreftrois#1{\definexref{¤#1}{{\the\titreun}.{\the\titredeux}.{\the\titretrois}}{}}
\else
\def\drefun#1{\definexref{prg#1}{{\the\titreun}}{}} 
\def\drefdeux#1{\definexref{prg#1}{{\the\titreun}.{\the\titredeux}}{}}
\def\dreftrois#1{\definexref{prg#1}{{\the\titreun}.{\the\titredeux}.{\the\titretrois}}{}}
\fi

%


  \long\def\propdeux#1#2#3#4{%
       \avance{\enonce}
       \leavevmode\edef\temp{#2}%
         \ifx\temp\empty 
          \else
           \definexref{#2}{#1~{\the\titreun.\the\enonce}}{enonces}
            \definexref{s#2}{{\the\titreun.\the\enonce}}{enonces}
             \fi
\smallskip
      \noindent{\bf#1\ {\bf\the\titreun.\the\enonce{#3}.}\enspace}{\sl#4\par}%
      \ifdim\lastskip<\medskipamount \removelastskip\penalty55\par \fi
   }

  \long\def\propun#1#2#3#4{%
      \avance{\enonce}
       \leavevmode\edef\temp{#2}%
        \ifx\temp\empty 
          \else
           \definexref{#2}{#1~{\the\enonce}}{enonces}
            \definexref{{s#2}}{{\the\enonce}}{enonces}
             \fi
   \par 
     \noindent{\bf#1\ {\bf\the\enonce{#3}.}\enspace}{\sl#4\par}%
     \ifdim\lastskip<\medskipamount \removelastskip\penalty55\medskip\fi
  }
  
  \long\def\prop#1#2#3#4{\ifnum\optionparag=1
                          \propdeux{#1}{#2}{\textfont1=\elevenib#3}{#4} \else\propun{#1}{#2}{\textfont1=\elevenib#3}{#4}\fi}

  \long\def\propt#1#2#3{\ifx\tpf\oui \prop{Th\'eo\-r\`eme}{#1}{#2}{#3}\par
                       \else\prop{Theorem}{#1}{#2}{#3}\par\fi}
  \long\def\Propt#1#2{\propt{#1}{}{#2}}
  \long\def\propl#1#2#3{\ifx\tpf\oui\prop{Lem\-me}{#1}{#2}{#3}\par
                         \else\prop{Lemma}{#1}{#2}{#3}\par\fi}
  \long\def\Propl#1#2{\propl{#1}{}{#2}}
  \long\def\propc#1#2#3{\ifx\tpf\oui\prop{Corol\-laire}{#1}{#2}{#3}\par
                         \else\prop{Corollary}{#1}{#2}{#3}\par\fi}

  \long\def\propd#1#2#3{\ifx\tpf\oui\prop{D\'efi\-nition}{#1}{#2}{#3}\par
                       \else\prop{Definition}{#1}{#2}{#3}\par\fi} 
  
  \long\def\proptd#1#2#3{\ifx\tpf\oui\prop{Th\'eor\`eme et d\'efi\-nition}{#1}{#2}{#3}\par
                       \else\prop{Theorem and definition}{#1}{#2}{#3}\par\fi}


  
  \newcount\optionparag\optionparag=1
  
  \long\def\section#1#2{\ifnum\optionparag=1 \paraun{#2} 
                        \else\goodbreak{\fontetitreun
  	                \Indentation{#1.\ }{#2}}\nobreak\smallskip\fi}

  \def\eqconstruct#1{\ifnum\optionparag=1{\the\titreun\hbox{$\cdot$}#1}\else{#1}\fi}

  
  
  \def\numref{oui}  
  
  \newcount\mesref\mesref=0 
  \def\defbib#1{\ifx\numref\oui\global\advance\mesref by 1 \definexref{#1}{{\the
                 \mesref}}{}\else\definexref{#1}{#1}{}\fi}
  \def\bibtem#1{\defbib{#1}\item{\citer{#1}}}
  \def\citer#1{[\ref{#1}]}
  \def\citeplus#1#2{[\ref{#1}; #2]}

  
  \font\seventeenmsa=msam10 at 17pt    
  \font\fourteenmsa=msam10 at 14pt
  \font\twelvemsa=msam10 at 12pt
  \font\tenmsa=msam10                 
  \font\ninemsa=msam10 at 9pt 
  \font\eightmsa=msam10 at 8pt 
  \font\sevenmsa=msam7 
  \font\sixmsa=msam10 at 6pt
  \font\fivemsa=msam5
  \newfam\msafam\textfont\msafam=\tenmsa\scriptfont\msafam=\sevenmsa\scriptscriptfont\msafam=\fivemsa
  
  \font\seventeenbb=msbm10 at 17pt     
  \font\fourteenbb=msbm10 at 14pt
  \font\twelvebb=msbm10 at 12pt
  \font\tenbb=msbm10                   
  \font\ninebb=msbm10 at 9pt 
  \font\eightbb=msbm10 at 8pt 
  \font\sevenbb=msbm7 
  \font\sixbb=msbm10 at 6pt
  \font\fivebb=msbm5 
  \newfam\bbfam\textfont\bbfam=\tenbb\scriptfont\bbfam=\sevenbb\scriptscriptfont\bbfam=\fivebb
  \def\bb{\fam\bbfam\tenbb}%

  \font\seventeenscaln=eusm10 at 17pt   
  \font\twelvescaln=eusm10 at 12pt
  \font\tenscaln=eusm10                
  \font\ninescaln=eusm10 scaled 900
  \font\eightscaln=eusm10 scaled 800
  \font\sevenscaln=eusm10 scaled 700
  \font\sixscaln=eusm10 scaled 600
   
  \newfam\scalnfam\textfont\scalnfam=\tenscaln\scriptfont\scalnfam=\sevenscaln\scriptscriptfont\scalnfam=\sixscaln
  \def\scaln{\fam\scalnfam\tenscaln}%
  \def\scal{\scaln}
  
  \font\tenscalb=eusb10                

  \font\sevenscalb=eusb10 scaled 700

  \newfam\scalbfam\textfont\scalbfam=\tenscalb\scriptfont\scalbfam=\sevenscalb
  %
  
  %
  %
  \font\fourteenrm=cmr12 scaled 1200
  \font\elevenrm=cmr10 at 11pt
  \font\twelverm=cmr12
  \font\ninerm=cmr9
  \font\eightrm=cmr8      
  \font\sevenrm=cmr7
  \font\sixrm=cmr6

  \font\seventeenpcap=cmcsc10 at 17pt
  \font\tenpcap=cmcsc10                        
  \font\ninepcap=cmcsc9
  \font\eightpcap=cmcsc8
  \font\sevenpcap=cmcsc10 scaled 700
  
  \newfam\pcapfam\textfont\pcapfam=\tenpcap\scriptfont\pcapfam=\sevenpcap
  \def\pcap{\fam\pcapfam\tenpcap}
  
  \font\seventeenrm=cmbx12 scaled 1400

  \font\fourteenbf=cmbx10 scaled 1400
  
  \font\twelvebf=cmbx12
  \font\elevenbf=cmbx10 at 11pt
  \font\ninebf=cmbx9  
  \font\eightbf=cmbx8
  \font\sixbf=cmbx6
  
  \font\tengot=eufm10                           
   
  \font\eightgot=eufm10 at 8truept 
  \font\sevengot=eufm7 
  \font\sixgot=eufm10 at 6 truept 
   
  \newfam\gotfam
  \textfont\gotfam=\tengot\scriptfont\gotfam=\sevengot\scriptscriptfont\gotfam=\sixgot
  \def\got{\fam\gotfam\tengot}%

  
  \def\tit{%
  \textfont0=\seventeenrm\scriptfont0=\tenrm\def\rm{\fam0\seventeenrm}%
  \textfont1=\seventeenib\scriptfont1=\twelveib%
  \textfont2=\seventeensy\scriptfont2=\twelvesy\scriptscriptfont2=\ninesy
  \textfont3=\seventeenex
  \textfont\itfam=\seventeenti
  \def\it{\fam\itfam\seventeenti}%
  \textfont\bbfam=\seventeenbb \scriptfont\bbfam=\twelvebb
  \def\bb{\fam\bbfam\seventeenbb}%
  \textfont\msafam=\seventeenmsa\scriptfont\msafam=\twelvemsa
  \textfont\scalnfam=\seventeenscaln
  \def\pcap{\fam\pcapfam\seventeenpcap}
  \normalbaselineskip=25pt\normalbaselines\rm}

  \font\seventeenti=cmbxti10 scaled 1680
  
  \font\fourteenti=cmbxti10 at 14pt
  
  \font\twelveti=cmbxti10 scaled 1200
  \font\eleventi=cmbxti10 at 11pt

  %
  %
  \font\twelveit=cmti12	
  \font\elevenit=cmti10 scaled 1100
  \font\nineit=cmti9
  \font\eightit=cmti8
  \font\sevenit=cmti7

  %
  %
 
 \font\seventeenib=cmmib10 scaled 1680
  \font\fourteenib=cmmib10 scaled 1400
  \font\twelveib=cmmib10 scaled 1200
  \font\elevenib=cmmib10 scaled 1100
  \font\tenib=cmmib10
\font\eightib=cmmib10 scaled 800
  \font\nineib=cmmib10 scaled 900
\font\sevenib=cmmib10 scaled 700
\font\sixib=cmmib10 scaled 600
\font\fiveib=cmmib10 scaled 500

\ifx\ITAN\oui
\else
\innernewfam\cmmibfam
\textfont\cmmibfam=\tenib
\scriptfont\cmmibfam=\sevenib
\scriptscriptfont\cmmibfam=\fiveib
\def\ib{\fam\cmmibfam\tenib}
\fi

  %
  %
  \font\twelvei=cmmi10 scaled 1200
  \font\eleveni=cmmi10 scaled 1100
  \font\ninei=cmmi9
  \font\eighti=cmmi8 
  \font\seveni=cmmi7 			                
  \font\sixi=cmmi6
  
  \font\ninesl=cmsl9                    
  \font\eightsl=cmsl8 
  \font\sevensl=cmsl10 at 7pt

  \font\ninett=cmtt9                    
  \font\eighttt=cmtt8
  \font\seventt=cmtt10 scaled 700

  \font\seventeensy=cmsy10 scaled 1680    
  \font\fourteensy=cmsy10 scaled 1400
  \font\twelvesy=cmsy10 scaled 1176
  
  \font\ninesy=cmsy9                      
  \font\eightsy=cmsy8
  \font\sixsy=cmsy6
  \font\seventeenex=cmex10 at 17pt
  \font\fourteenex=cmex10 at 14pt
  \font\twelveex=cmex10 at 12pt
  \font\nineex=cmex10 at 9pt
  \font\eightex=cmex10 at 8pt
  \font\sevenex=cmex10 at 7pt
  \font\sixex=cmex10 at 6pt
  \font\fiveex=cmex10 at 5pt
  
   
  \font\fourteengp=cmmi10 at 14pt
  
  \font\twelvegp=cmmib10 at 12pt
  \font\elevengp=cmmib10 at 11pt
  \font\tengp=cmmib10                          
  \font\ninegp=cmmib10 at 9pt 
  \font\eightgp=cmmib8 
   
  \font\sixgp=cmmib6


  \def\gponze{\textfont0=\elevenbf\scriptfont0=\eightbf\scriptscriptfont0=\sixbf
           \textfont1=\elevengp\scriptfont1=\eightgp\scriptscriptfont1=\sixgp}
  \def\gpdouze{\textfont0=\twelvebf\scriptfont0=\tenbf\scriptscriptfont0=\ninebf
           \textfont1=\twelvegp\scriptfont1=\tengp\scriptscriptfont1=\ninegp}        
  
 \def\gpquatorze{\textfont0=\fourteenbf\scriptfont0=\twelvebf\scriptscriptfont0=\elevenbf
           \textfont1=\fourteengp\scriptfont1=\twelvegp\scriptscriptfont1=\elevengp}

  
  \expandafter\chardef\csname pre amssym.def at\endcsname=\the\catcode`\@
  \catcode`\@=11
  \def\undefine#1{\let#1\undefined}
  \def\newsymbol#1#2#3#4#5{\let\next@\relax
   \ifnum#2=\@ne\let\next@\msafam@\else
   \ifnum#2=\tw@\let\next@\bbfam@\fi\fi
   \mathchardef#1="#3\next@#4#5}
  \def\mathhexbox@#1#2#3{\relax
   \ifmmode\mathpalette{}{\m@th\mathchar"#1#2#3}%
   \else\leavevmode\hbox{$\m@th\mathchar"#1#2#3$}\fi}
  \def\hexnumber@#1{\ifcase#1 0\or 1\or 2\or 3\or 4\or 5\or 6\or 7\or 8\or
   9\or A\or B\or C\or D\or E\or F\fi}
  
  \def\setboxz@h{\setbox\z@\hbox}
  \def\wdz@{\wd\z@}
  \def\boxz@{\box\z@}
  
  \edef\msafam@{\hexnumber@\msafam}
  \mathchardef\dabar@"0\msafam@39
  
  \edef\bbfam@{\hexnumber@\bbfam}
  \def\widehat#1{\setboxz@h{$\m@th#1$}%
   \ifdim\wdz@>\tw@ em\mathaccent"0\bbfam@5B{#1}%
   \else\mathaccent"0362{#1}\fi}
  \def\widetilde#1{\setboxz@h{$\m@th#1$}%
   \ifdim\wdz@>\tw@ em\mathaccent"0\bbfam@5D{#1}%
   \else\mathaccent"0365{#1}\fi}
  \newsymbol\leqq 1335          
  \newsymbol\leqslant 1336
  \newsymbol\lessgtr 1337       
  \newsymbol\backprime 1038     
  \newsymbol\risingdotseq 133A  
  \newsymbol\fallingdotseq 133B 
  \newsymbol\succcurlyeq 133C   
  \newsymbol\geqq 133D          
  \newsymbol\geqslant 133E
  \newsymbol\nmid 232D
  \newsymbol\nexists 2040
  \newsymbol\smallsetminus 2272
  \newsymbol\varnothing 203F 
  \catcode`\@=\active

  \catcode`\@=11
  \newcount\typofr\typofr=1
  
  \catcode`\;=\active
  \def;{\ifnum\typofr=1\relax\ifhmode\ifdim\lastskip>\z@\unskip\fi
     \kern.2em\fi\string;\else\string;\fi}
  
  \catcode`\:=\active
  \def:{\ifnum\typofr=1\relax\ifhmode\ifdim\lastskip>\z@\unskip\fi
  \penalty\@M\ \fi\string:\else\string:\fi}
  
  \catcode`\!=\active
  \def!{\ifnum\typofr=1\relax\ifhmode\ifdim\lastskip>\z@\unskip\fi
     \kern.2em\fi\string!\else\string!\fi}
  
  \catcode`\?=\active
  \def?{\ifnum\typofr=1\relax\ifhmode\ifdim\lastskip>\z@\unskip\fi
     \kern.2em\fi\string?\else\string?\fi}

  \def\francais{\typofr=1\def\tpf{oui}}
  \def\anglais{\typofr=2\def\tpf{non}\def\english{oui}}
  \def\oui{oui}
  \francais
  
  \catcode`\@=12
  

%
\def\raye #1|{\leavevmode\setbox1=\hbox{#1}%
\raise .5pt\hbox to \wd1{\xleaders\hbox{\rge{$ \sct / $}%
\kern 1pt}\hfill\kern -1pt }\kern -\wd1 \unhbox1\relax }

\def\barre #1|{\leavevmode\setbox1=\hbox{#1}%
\rlap{\color{red}\vrule height 2.4pt depth -1.2pt width \wd1}\color{black} \unhbox1\relax}

  

  
  \def\og{\leavevmode\raise.24ex\hbox{$\scriptscriptstyle\langle\!\langle\>$}}    
  \def\fg{\leavevmode\raise.24ex\hbox{$\scriptscriptstyle\>\rangle\!\rangle$}}    

  \def\d{\,{\rm d}}

  \def\r{{\bb R}}
  
  \def\N{{\bb N}}

  \def\K{{\scal K}}

  \def\M{{\scal M}}
  
  \def\O{{\scal O}}
  \def\P{{\scaln P}}

  \def\frac#1#2{{#1\over #2}}
  \def\di#1#2{\sct#1\atop{\sct#2}}

  \def\numero{n$^{\rm o}\thinspace$}

  \def\qedbox{$\rlap{$\sqcap$}\sqcup$}           
  \def\qed{\nobreak\hfill\penalty250 \hbox{}\nobreak\hfill\qedbox\par }

  \def\numero{n$^{\rm o}\thinspace$}

  \def\pnu{p^\nu}

  \def\¤{\S\thinspace}

  \def\¥{$\bullet$ }
  
  
  \def\e{{\rm e}}

  \def\epsilon{\varepsilon}

  \def\phi{\varphi}
  \def\theta{\vartheta}
  \def\rho{\varrho}
  \def\dm{{\textstyle{1\over 2}}}
  \def\txt{\textstyle}
  \def\dsp{\displaystyle}
  \def\sct{\scriptstyle}
  \def\pf{\noi{\it Proof. }}
  \def\nid{\ifnum\typofr=1\par\noindent{\it D\'emonstration. }\else\pf\fi}
  \def\noi{\noindent}
  \def\rem{\ifnum\typofr=1\noi{\it Remarque.}\ \else\noi{\it Remark.}\ \fi}
  \def\rems{\ifnum\typofr=1\noi{\it Remarques.}\ \else\noi{\it Remarks.}\ \fi}
  \def\re{{\Re e\,}}

  \def\sset{\smallsetminus}

  \def\1{{\bf 1}}
  \def\|{\Vert}

  \def\leq{\leqslant}
  \def\geq{\geqslant}

  \def\eg{{e.g.}}
  \newsymbol\subsetneqq 2324
  \newsymbol\subsetneq 2328

  \def\log{\mathop{\rm log}\nolimits}
  \def\ft#1#2{{\txt{#1\over #2}}}



\def\Vbs#1{\bigg|#1\bigg|}


  \def\pmb#1{\setbox0=\hbox{#1}%
  \kern-.025em\copy0\kern-\wd0\kern.05em\copy0\kern-\wd0\kern-.025em\raise .0433em\box0 }

  
  \skewchar\eighti='177 \skewchar\sixi='177
  \skewchar\eightsy='60 \skewchar\sixsy='60
  
  \def\eightpoint{%
  \textfont0=\eightrm\scriptfont0=\sixrm\scriptscriptfont0=\fiverm
  \def\rm{\fam0\eightrm}%
  \textfont1=\eighti\scriptfont1=\sixi
  \scriptscriptfont1=\fivei\def\oldstyle{\fam1\seveni}%
  \textfont2=\eightsy\scriptfont2=\sixsy\scriptscriptfont2=\fivesy
  \textfont3=\eightex\scriptfont3=\sixex
  \textfont\itfam=\eightit
  \def\it{\fam\itfam\eightit}%
  \textfont\slfam=\eightsl
  \def\sl{\fam\slfam\eightsl}%
  \textfont\bbfam=\eightbb \scriptfont\bbfam=\sixbb\scriptscriptfont\bbfam=\fivebb
  \def\bb{\fam\bbfam\eightbb}%
  \textfont\msafam=\eightmsa\scriptfont\msafam=\sixmsa
  \textfont\scalnfam=\eightscaln
  \def\scaln{\fam\scalnfam\eightscaln}
  \textfont\ttfam=\eighttt
  \def\tt{\fam\ttfam\eighttt}%
\textfont\gotfam=\eightgot
  \textfont\bffam=\eightbf\scriptfont\bffam=\sixbf\scriptscriptfont\bffam=\fivebf
  \def\bf{\fam\bffam\eightbf}%
  \ifx\ITAN\oui\else\textfont\cmmibfam=\eightib
       \scriptfont\cmmibfam=\sixib
        \scriptscriptfont\cmmibfam=\fiveib
         \def\ib{\fam\cmmibfam\eightib}
   \fi
  \textfont\pcapfam=\eightpcap
  \def\pcap{\fam\pcapfam\eightpcap}
  \abovedisplayskip=2pt plus2pt minus 2pt
  \belowdisplayskip=2pt plus1pt minus 2pt
  \abovedisplayshortskip= 1pt plus 2pt minus 1pt
  \belowdisplayshortskip= 1pt plus 2pt minus 1pt
  \smallskipamount=2pt plus 1pt minus 2pt
  \medskipamount=3pt plus 2pt minus 2pt
  \bigskipamount=7pt plus 3pt minus 3pt
  \setbox\strutbox=\hbox{\vrule height 5pt depth 2pt width 0pt}%
  \normalbaselineskip=9pt\normalbaselines\rm}

  \def\({\left(}
  \def\){\right)}
  
  \def\footnoterule{\kern -2pt\hrule width 7truecm\kern 2.4pt}
  
  \def\xnotedef#1{\definexref{#1}{\noexpand\number\footnotenumber}{Note}}%

  
  
  \def\ninepoint{%
  \textfont0=\ninerm\scriptfont0=\sixrm\scriptscriptfont0=\fiverm
  \def\rm{\fam0\ninerm}%
  \textfont1=\ninei\scriptfont1=\sixi
  \scriptscriptfont1=\fivei\def\oldstyle{\fam1\ninei}%
  \textfont2=\ninesy\scriptfont2=\sixsy\scriptscriptfont2=\fivesy
  \textfont3=\nineex\scriptfont3=\sixex
  \textfont\itfam=\nineit
  \def\it{\fam\itfam\nineit}%
  \textfont\slfam=\ninesl
  \def\sl{\fam\slfam\ninesl}%
  \textfont\bbfam=\ninebb\scriptfont\bbfam=\sixbb\scriptscriptfont\bbfam=\fivebb
  \def\bb{\fam\bbfam\ninebb}%
  \textfont\msafam=\ninemsa\scriptfont\msafam=\sixmsa\scriptscriptfont\msafam=\fivemsa
  \textfont\scalnfam=\ninescaln
  \def\scaln{\fam\scalnfam\ninescaln}
  \textfont\ttfam=\ninett
  \def\tt{\fam\ttfam\ninett}%
  \textfont\bffam=\ninebf\scriptfont\bffam=\sixbf\scriptscriptfont\bffam=\fivebf
  \def\bf{\fam\bffam\ninebf}%
  \abovedisplayskip=3pt plus2pt minus 2pt
  \belowdisplayskip=3pt plus1pt minus 2pt
  \abovedisplayshortskip= 2pt plus 2pt minus 1pt
  \belowdisplayshortskip= 2pt plus 2pt minus 1pt
  \smallskipamount=2pt plus 1pt minus 2pt
  \medskipamount=3pt plus 2pt minus 2pt
  \bigskipamount=7pt plus 3pt minus 3pt
  \setbox\strutbox=\hbox{\vrule height 5pt depth 2pt width 0pt}%
  \normalbaselineskip=11pt plus.3pt minus.2pt\normalbaselines\rm}

  \def\sevenpoint{%
  \textfont0=\sevenrm\scriptfont0=\sixrm\scriptscriptfont0=\fiverm
  \def\rm{\fam0\sevenrm}%
  \textfont1=\seveni\scriptfont1=\sixi
  \scriptscriptfont1=\fivei\def\oldstyle{\fam1\seveni}%
  \textfont2=\sevensy\scriptfont2=\sixsy\scriptscriptfont2=\fivesy
  \textfont3=\sevenex\scriptfont3=\fiveex
  \textfont\itfam=\sevenit
  \def\it{\fam\itfam\sevenit}%
  \textfont\slfam=\sevensl
  \def\sl{\fam\slfam\sevensl}%
  \textfont\bbfam=\sevenbb \scriptfont\bbfam=\sixbb\scriptscriptfont\bbfam=\fivebb
  \def\bb{\fam\bbfam\sevenbb}%
  \textfont\msafam=\sevenmsa\scriptfont\msafam=\sixmsa
  \textfont\scalnfam=\sevenscaln
  \def\scaln{\fam\scalnfam\sevenscaln}
  \textfont\bffam=\sevenbf\scriptfont\bffam=\sixbf\scriptscriptfont\bffam=\fivebf
  \def\bf{\fam\bffam\sevenbf}%
  \textfont\ttfam=\seventt
  \abovedisplayskip=2pt plus2pt minus 2pt
  \belowdisplayskip=2pt plus1pt minus 2pt
  \abovedisplayshortskip= 1pt plus 2pt minus 1pt
  \belowdisplayshortskip= 1pt plus 2pt minus 1pt
  \smallskipamount=2pt plus 1pt minus 2pt
  \medskipamount=3pt plus 2pt minus 2pt
  \bigskipamount=7pt plus 3pt minus 3pt
  \setbox\strutbox=\hbox{\vrule height 5pt depth 2pt width 0pt}%
  \normalbaselineskip=9pt\normalbaselines\rm}

 \def\onzepts{%
 \textfont0=\elevenrm\scriptfont0=\ninerm
 \textfont1=\eleveni\scriptfont1=\ninei
}

\def\douzepts{%
  \textfont0=\twelverm\scriptfont0=\tenrm\def\rm{\fam0\twelverm}%
  \textfont1=\twelvei\scriptfont1=\teni%
  \textfont2=\twelvesy\scriptfont2=\tensy\scriptscriptfont2=\eightsy
  \textfont3=\twelveex
  \textfont\itfam=\twelveti
  \def\it{\fam\itfam\twelveti}%
  \textfont\bffam=\twelvebf\scriptfont\bffam=\tenbf\scriptscriptfont\bffam=\eightbf
  \def\bf{\fam\bffam\twelvebf}%
  \textfont\bbfam=\twelvebb \scriptfont\bbfam=\tenbb
  \def\bb{\fam\bbfam\twelvebb}%
  \textfont\msafam=\twelvemsa\scriptfont\msafam=\tenmsa
  \textfont\scalnfam=\twelvescaln
  \normalbaselineskip=15pt\normalbaselines\rm}

\def\quatorzepts{%
  \textfont0=\fourteenrm\scriptfont0=\twelverm\def\rm{\fam0\fourteenrm}%
  \textfont1=\fourteenib\scriptfont1=\twelveib%
  \textfont2=\fourteensy\scriptfont2=\twelvesy\scriptscriptfont2=\tensy
  \textfont3=\fourteenex
  \textfont\itfam=\fourteenti
  \def\it{\fam\itfam\fourteenti}%
  \textfont\bffam=\fourteenbf\scriptfont\bffam=\twelvebf\scriptscriptfont\bffam=\tenbf
  \def\bf{\fam\bffam\fourteenbf}%
  \textfont\bbfam=\fourteenbb \scriptfont\bbfam=\twelvebb
  \def\bb{\fam\bbfam\fourteenbb}%
  \textfont\msafam=\fourteenmsa\scriptfont\msafam=\twelvemsa
  \textfont\scalnfam=\twelvescaln
  \normalbaselineskip=18pt\normalbaselines\rm}


\def\AA{{\it Acta Arith.}}

\def\picture #1 by #2 (#3){\leavevmode\vbox to #2{
     \hrule width #1 height 0pt depth 0pt
      \vfill
       \special{picture #3}}}

\def\illustration #1 by #2 (#3) scaled #4{\dimen1=#2
  \divide\dimen1 by 1000
  \multiply\dimen1 by #4
  \vtop to \dimen1{\dimen1=#1
  \divide\dimen1 by 1000
  \multiply\dimen1 by #4
  \hsize=\dimen1\vss
  \noindent\special{illustration #3 scaled #4}}}

\ifx\couleurs\oui

\fi

\ifx\optionkeymacros\undefined\else \fi

\catcode`\Œ=\active\defŒ{{\aa}}       
\catcode`\º=\active\defº{\int}        
\catcode`\=\active\def{\c c}        
\catcode`\¶=\active\def¶{\partial}    
\catcode`\Ä=\active\defÄ{\oint}       
\catcode`\Æ=\active\defÆ{\triangle}   
\catcode`\Â=\active\defÂ{\neg}        
\catcode`\µ=\active\defµ{\mu}         
\catcode`\¿=\active\def¿{{\o}}        
\catcode`\¹=\active\def¹{\pi}         
\catcode`\Ï=\active\defÏ{{\oe}}       
\catcode`\§=\active\def§{{\ss}}       
\catcode`\ =\active\def {\dagger}     
\catcode`\Ã=\active\defÃ{\sqrt}       
\catcode`\·=\active\def·{\Sigma}      
\catcode`\Å=\active\defÅ{\approx}     
\catcode`\½=\active\def½{\Omega}      
\catcode`\£=\active\def£{{\it\$}}     
\catcode`\°=\active\def°{\infty}      
\catcode`\¤=\active\def¤{{\S}}        
\catcode`\¦=\active\def¦{{\P}}        
\catcode`\¥=\active\def¥{\bullet}     
\catcode`\»=\active\def»{\leavevmode\raise.585ex\hbox{\b a}}      
\catcode`\¼=\active\def¼{\leavevmode\raise.6ex\hbox{\b o}}        
\catcode`\­=\active\def­{\not=}       
\catcode`\²=\active\def²{\leq}        
\catcode`\³=\active\def³{\geq}        
\catcode`\Ö=\active\defÖ{\div}        
\catcode`\É=\active\defÉ{{\dots}}     
\catcode`\¾=\active\def¾{{\ae}}       
\catcode`\Ç=\active\defÇ{\og}         
\catcode`\Ò=\active\defÒ{``}          
\catcode`\Á=\active\defÁ{!`}          
\catcode`\¢=\active\def¢{\rlap/c}     
\catcode`\Ô=\active\defÔ{`}           
\catcode`\Õ=\active\defÕ{'}           


\catcode`\=\active\def{{\AA}}       
\catcode`\'=\active\def'{\c C}        
\catcode`\¯=\active\def¯{{\O}}        
\catcode`\¸=\active\def¸{\Pi}         
\catcode`\Î=\active\defÎ{{\OE}}       
\catcode`\®=\active\def®{{\AE}}       
\catcode`\×=\active\def×{\diamond}    
\catcode`\¡=\active\def¡{\accent'27}  
\catcode`\Ó=\active\defÓ{''}          
\catcode`\±=\active\def±{\pm}         
\catcode`\È=\active\defÈ{\fg}         
\catcode`\À=\active\defÀ{?`}          
\catcode`\Ð=\active\defÐ{--}          
\catcode`\Ñ=\active\defÑ{---}         


\catcode`\Š=\active\defŠ{\"a}        
\catcode`\'=\active\def'{\"e}        
\catcode`\•=\active\def•{\"{\i}}     
\catcode`\š=\active\defš{\"o}        
\catcode`\Ÿ=\active\defŸ{\"u}        
\catcode`\Ø=\active\defØ{\"y}        
\catcode`\å=\active\defå{\^A}        
\catcode`\€=\active\def€{\"A}        
\catcode`\…=\active\def…{\"O}        
\catcode`\†=\active\def†{\"U}        
\catcode`\‡=\active\def‡{\'a}        
\catcode`\Ž=\active\defŽ{\'e}        
\catcode`\'=\active\def'{\'{\i}}     
\catcode`\—=\active\def—{\'o}        
\catcode`\œ=\active\defœ{\'u}        
\catcode`\ƒ=\active\defƒ{\'E}        
\catcode`\æ=\active\defæ{\^E}        
\catcode`\é=\active\defé{\`E}        %
\catcode`\ˆ=\active\defˆ{\`a}        
\catcode`\=\active\def{\`e}        
\catcode`\"=\active\def"{\`{\i}}     
\catcode`\˜=\active\def˜{\`o}        
\catcode`\=\active\def{\`u}        
\catcode`\Ë=\active\defË{\`A}        
\catcode`\‹=\active\def‹{\~a}        
\catcode`\–=\active\def–{\~n}        
\catcode`\›=\active\def›{\~o}        
\catcode`\Ì=\active\defÌ{\~A}        
\catcode`\"=\active\def"{\~N}        
\catcode`\Í=\active\defÍ{\~O}        
\catcode`\‰=\active\def‰{\^a}        
\catcode`\=\active\def{\^e}        
\catcode`\"=\active\def"{\^{\i}}     
\catcode`\™=\active\def™{\^o}        
\catcode`\ž=\active\defž{\^u}        

\let\optionkeymacros\null

\vsize=254truemm
\voffset=-5truemm
\dateheure
\anglais
\scriptscriptfont0=\fiverm

    \font\tenrsfs=rsfs7 at 10pt

    \font\sevenrsfs=rsfs7
    \font\sixrsfs=rsfs7 at 6pt

\newfam\rsfsfam\textfont\rsfsfam=\tenrsfs\scriptfont\rsfsfam=\sevenrsfs\scriptscriptfont\rsfsfam=\sixrsfs

\def\ga{{\got a}}
\def\gb{{\got b}}
\def\gc{{\got c}}

\def\gh{{\got h}}

\def\gpp{{\got p}}

\def\fmus{mul\-ti\-pli\-cative func\-tions}

\hautspages{G. Tenenbaum}{On effective mean-values of arithmetic functions}
\vskip-3mm
\titrecentre{On effective mean-values of arithmetic functions}
\bigskip
\centerline{GŽrald Tenenbaum}
\smallskip
{\leftskip110mm
\obeylines
\it To Krishna Alladi,
inspired founder of
The Ramanujan Journal
\par }
\bigskip
{\eightpoint\leftskip1cm\rightskip1cm
\noi{\bf Abstract.} Let $r,\,f$ be multiplicative functions with $r\geqslant 0$, $f$ is complex valued, $|f|\leqslant r$, and $r$ satisfies some standard growth hypotheses. Let $x$ be large, and assume that, for some real number~$\tau$, the quantities $r(p)-\re\{f(p)/p^{i\tau}\}$ are small in various appropriate average senses over the set of prime numbers not exceeding $x$. We derive from recent effective mean-value estimates an effective comparison theorem between the mean-values of $f$ and of $r$ on the set of integers~$\leqslant x$. We also provide effective estimates for certain weighted moments of additive functions and for sifted mean-values of non-negative multiplicative functions.
 \PAR
\medskip
\noi
{\bf Keywords:} Quantitative estimates, comparison theorems, multiplicative functions, effective mean-value estimates, Wirsing's theorem, H‡lasz' theorem, weighted distribution of additive functions, weighted moments of additive functions..\par
\medskip
\noi \bf 2020 Mathematics Subject Classification: \rm primary   11K65, 11N37; secondary 11N25, 11N60, 11N64.\par }
\bigskip
\medskip
\paraun{Introduction}
There is an abundant literature on estimates for mean-values multiplicative functions, usually appearing in the form of summatory functions
$$M(x;f):=\sum_{n\leqslant x}f(n)\qquad (x\geqslant 1).$$
One of the most useful class of results in this topic is that of {\it comparison theorems}, evaluating ratios $M(x;f)/M(x;r)$ where $r$ is a majorant of $|f|$. The first general theorems of this type are due to Wirsing \citer{Wi67} and Hal‡sz \citer{Ha68}, and we refer to \citer{Te17} for a more complete account of the literature.
\par 
This note is devoted to describing some consequences of  effective estimates of the above kind obtained in \citer{Te17} and which we now state. 
\par 
Let $\M(A,B)$ designate the class of those complex multiplicative functions $f$ satisfying 
$$\max_{p}|f(p)|\leqslant A,\quad  \sum_{p,\,\nu\geqslant 2}{|f(p^\nu)|\log p^\nu\over p^\nu}\leqslant B,\eqdef{C0}
$$
where, here and in the sequel, the letter $p$ denotes a prime number. 
Define furthermore $w_f:=1$ if $f$ is real, $w_f:=\dm$ if $f$ assumes some non real values, and write
$$Z(x;f):=\sum_{p\leqslant x}{f(p)\over p}\cdot$$
\par 
The following two statements are established in \citer{Te17}. Here and throughout the notation $]u,v]$ indicates a semi-open interval, including $v$ but not $u$. 
\propt{cascv}{ (\citer{Te17})}{Let $$\eqalign{&\ga\in]0,\dm],\ \gb\in\,[\ga,1[,\  A\geqslant2\gb,\  B>0,\,
x\geqslant \e^4,\, \varepsilon=\varepsilon_x\in\big]1/\sqrt{\log x}, \dm\big],\cr &\varrho=\varrho_x\in\,[2\gb,A],\
\gpp:={\pi\varrho\over A},\ \beta:=1-{\sin\gpp\over \gpp},\quad\gh:={1-\gb\over \min(1,\varrho)-\gb}\cdot\cr}\eqdef{defpar}$$  
Assume that the \fmus\    $f$, $r$, satisfy $r\in\M(A,B)$, $|f|\leqslant r$, and
$$\leqalignno{\sum_{p\leqslant x}{r(p)-\re f(p)\over p}&\leqslant \ft12\beta\gb\log (1/\epsilon), &\eqdef{pqr}\cr
\sum_{x^\varepsilon<p\leqslant y}{\{r(p)-\re f(p)\}^\gh\log p\over p}&\ll \varepsilon^{\delta_1\gh}\log y\qquad (x^\varepsilon<y\leqslant x), &\eqdef{cdpr1}\cr
\sum_{p\leqslant y}{\{r(p)-\varrho\}\log p\over p}&\ll\varepsilon\log y\qquad 
(x^{\epsilon}< y\leqslant x),&\eqdef{C2}\cr
}$$ 
where
$\delta_1\in]0,\ft23\beta\gb]$. Then we have 
$$M(x;f)={\e^{-\gamma\varrho}x\over \Gamma(\varrho)\log
x}\bigg\{\prod_{p}\sum_{\pnu\leqslant 
x} {f(p^\nu)\over p^\nu}+O\Big(\varepsilon^\delta\,\e^{Z(x;f)}\Big)\bigg\},\eqdef{moy}$$
with $
\delta:=w_f\delta_1$, and where
the implicit constant in \eqref{moy} depends at most on $A$, $B$, $\ga$, 
 $\gb$, and on the implicit constants in \eqref{cdpr1} and \eqref{C2}. }
\medskip
 As noted in \citer{Te17}, hypothesis \eqref{cdpr1} is trivially  implied by the  condition
$$\sum_{x^\varepsilon<p\leqslant x}{\{r(p)-\re f(p)\}^\gh\over p}\ll  \epsilon^{\delta_1\gh},\eqdef{cdpr}$$
and, of course, also by the uniform bound 
$$\max_{x^\varepsilon<p\leqslant x}\{r(p)-\re f(p)\}\ll \varepsilon^{\delta_1}.\eqdef{cdpr2}$$
\par 
We also quote from \citer{Te17} the remark that, while the assumptions of \ref{cascv} imply
$$\prod_{p}\sum_{\pnu\leqslant 
x} {f(p^\nu)\over p^\nu}\ll\e^{Z(x;f)},$$ 
the two sides of \eqref{moy} have the same order of magnitude if
$$\min_{p,x}\Vbs{\sum_{0\leqslant \nu\leqslant (\log x)/\log p} {f(p^\nu)\over p^\nu}}\gg1.\eqdef{hypgen}$$
Under this condition generically satisfied, formula \eqref{moy} becomes
$$M(x;f)=\big\{1+O(\varepsilon^\delta)\big\}{\e^{-\gamma\varrho}x\over \Gamma(\varrho)\log
x}\prod_{p}\sum_{\pnu\leqslant 
x} {f(p^\nu)\over p^\nu}\cdot\eqdef{moy2}$$
\par 
For the next statement, we introduce the notation
$$\beta_0=\beta_0(\gb,A):=1-{\sin(2\pi\gb/A)\over 2\pi\gb/A},\qquad  \delta_0(\gb)=\delta_0(\gb,A):=\ft13\gb\beta_0.\eqdef{delta0}$$\par 
\vskip-5mm
 \propt{fr}{ (\citer{Te17})}{Let $$\ga\in]0,\ft14], \quad\gb\in[\ga,\ft12[,  \quad A\geqslant 2\gb,\quad B>0, \quad  x\geqslant \e^4, \quad 1/\sqrt{\log x}<\varepsilon\leqslant \dm.$$ Assume that the \fmus\ $f$, $r$, such that \hbox{$r\in\M(2A,B)$}, $|f|\leqslant r$,  satisfy conditions \eqref{pqr}, \eqref{cdpr1} with $\gh:=(1-\gb)/\gb$, \eqref{cdpr} with $\gh=1$, and 
$$\sum_{y<p\leqslant y^{1+\varepsilon_1}}{r(p)\log p\over p}\geqslant 4\gb\varepsilon_1\log y\qquad \big(\e^{1/\varepsilon_1}\leqslant y\leqslant x^{1/(1+\varepsilon_1)}\big)\eqdef{moyrp}$$
where $\varepsilon_1:=\sqrt{\varepsilon}$. Assume furthermore that $\delta_1\in]0,\delta_0(\gb)] $. Then we have
$$M(x;f)=M(x;r)\prod_{p}{\sum_{\pnu\leqslant x}f(\pnu)/\pnu\over \sum_{\pnu\leqslant x}r(\pnu)/\pnu}+O\bigg({x\,\varepsilon^{\delta}\e^{Z(x;r)-\gc Z(x;|f|-f)}\over \log x}\bigg)\eqdef{Mf/Mr}$$
where $\delta:=w_f\delta_1$, $\gc:=\gb/A$. Moreover, the above formula persists without requiring \eqref{cdpr} to hold with $\gh=1$ provided $\min_{x^\varepsilon\leqslant p\leqslant x}r(p)\geqslant 4\gb$.
The implicit constant in \eqref{Mf/Mr} depends at most on $A$, $B$, $\ga$, $\gb$, as well as on the implicit constants in \eqref{cdpr1} and if need be in \eqref{cdpr}.}
\par 
The main novelty in the statements of Theorems \ref{scascv} and \ref{sfr} is effectivity in hypothesis \eqref{cdpr1}, which finds its counterpart in the estimates for the remainder terms. The introduction of parameters $A$, $B$, $\ga$, $\gb$, is purely contingent and is designed to facilitate applicability.
\bigskip 
\bigskip
\paraun{An effective comparison theorem of Wirsing type}
\ref{fr} provides an estimate for $M(x;f)/M(x;r)$ when $r(p)-\re f(p)$ is small in suitable respects. However, under the assumption that, for suitable $\varrho>0$, we have $$\sum_{p\leqslant x}{\{r(p)-\varrho\}\log p\over p}=o\big(\log x\big),\eqdef{cW}$$ Wirsing's theorem \citer{Wi67} also provides, via partial summation,  an asymptotic formula for $M(x;f)$ whenever
the condition 
$$\sum_{p}{r(p)-\re\{f(p)/p^{i\tau}\}\over p}<\infty\eqdef{condcv}$$
holds for some real number $\tau$, necessarily unique. A further extension has recently been established under condition \eqref{condcv} by Indlekofer \& Kaya \citer{IK25},  assuming  moreover that $r=r_1*r_2$, $r_1\geqslant 0$, $r_2\geqslant 0$, with $r_1$~satisfying~\eqref{cW}, and the~$r_j$ are multiplicative functions with disjoint supports on the set of prime powers.
\par \goodbreak
While, both in the classical framework and in the setting considered in \citer{IK25}, the transition from $\tau=0$ to general $\tau$ is obtained through simple partial integration, the deduction is not straightforward when an effective estimate is aimed at. The following result, which is a consequence of \ref{fr}, furnishes the desired extension under the much weaker hypothesis \eqref{moyrp}. For given $\tau\in\r$, we write $f_\tau(n):=f(n)/n^{i\tau}$ $(n\geqslant 1)$ and recall notation \eqref{delta0} for $\delta_0(\gb)$.
\par \vskip-3mm
\Propt{thWg}{Let $$\ga\in]0,\ft14], \quad\gb\in[\ga,\ft12[,  \quad A\geqslant 2\gb,\quad B>0, \quad x\geqslant \e^4, \quad 1/\sqrt{\log x}<\varepsilon\leqslant \dm,\quad\tau\in\r.$$ Assume that the \fmus\ $f$, $r$, with \hbox{$r\in\M(2A,B)$}, $|f|\leqslant r$,  are such that conditions \eqref{pqr}, \eqref{cdpr1} with $\gh:=(1-\gb)/\gb$, and \eqref{cdpr} with $\gh=1$ are satisfied with $f_\tau$ in place of $f$. Suppose furthermore that \eqref{moyrp} holds and
that $\delta_1\in]0,\delta_0(\gb)] $. Then we have
$$M(x;f)={x^{i\tau}M(x;r)\over 1+i\tau}\prod_{p}{\sum_{\pnu\leqslant x}f(\pnu)/p^{\nu(1+i\tau)}\over \sum_{\pnu\leqslant x}r(\pnu)/\pnu}+O\Big(\varepsilon^\delta M(x;r)\Big),\eqdef{Mf/Mr-g}$$
where $\delta:=w_f\delta_1$.
Moreover, the above formula persists without requiring \eqref{cdpr} to hold with $\gh=1$ provided $\min_{x^\varepsilon\leqslant p\leqslant x}r(p)\geqslant 4\gb$.
The implicit constant in \eqref{Mf/Mr-g} depends at most on $A$, $B$, $\ga$, $\gb$, $\tau$, as well as on the implicit constants in \eqref{cdpr1} and if need be in \eqref{cdpr}.}
\medskip
For simplicity, we have omitted in the above statement the potential sharpening of the error term involving $Z(x;|f|-f)$ and appearing in \eqref{Mf/Mr}. With some extra care, it could be reinserted.
\medskip
We start with a lemma potentially useful in other contexts.
\Propl{Mzr}{Assume $r\in\M(2A,B)$, $r\geqslant 0$, and that \eqref{moyrp} holds with $\varepsilon$ sufficiently small. Then 
$$M(z;r)\asymp {z\e^{Z(z;r)}\over \log z}\qquad (x^{2\varepsilon_1}<z\leqslant x).\eqdef{Mxras}$$}
\nid
The corresponding upper bound is standard and follows from \citer{HR79} or \citer{Sh80}.
To establish the lower bound,  define $$\varepsilon_2:=\varepsilon\varepsilon_1=\varepsilon^{3/2},\quad\K:=\bigg[{\log (\varepsilon_2\log x)\over \log (1+\varepsilon_1)}, {\log_2x\over \log (1+\varepsilon_1)}-1\bigg]\cap\N,$$ and apply \eqref{moyrp} with $y=y_k:=\exp\{(1+\varepsilon_1)^k\}$ for $k\in\K$ to get
$$\sum_{y_k<p\leqslant y_{k+1}}{r(p)\log p\over p}=\gb_k\varepsilon_1\log y_k\qquad \big(k\in\K\big), \eqdef{moylocrp}$$ 
with $$4\gb\leqslant \gb_k\leqslant 2A+O\Big(\e^{-\sqrt{\log y_k}}\Big).$$
\par 
Then define $$s(p):=2\gb r(p)/\gb_k\quad \big(y_k<p\leqslant y_{k+1},\,k\in\K\big), \quad s(p):=0\quad\Big(p\in[2,x]\sset\cup_{k\in\K}\,]y_k,y_{k+1}]\Big).\eqdef{defsp}$$ Thus $0\leqslant s(p)\leqslant \dm r(p)$ for all $p\leqslant x$. By summation, it follows that
$$\sum_{p\leqslant y}{\{s(p)-2\gb\}\log p\over p}\ll\varepsilon_1\log y\qquad (x^{\varepsilon}<y\leqslant x).\eqdef{regs}$$  
Next, define two multiplicative functions $s_1$ and $t_1$, supported on squarefree integers, by the  formulae $s_1(p):=s(p)$, $t_1(p):=r(p)-s_1(p)$,  and put $r_1:=s_1*t_1$.
\par 
 For $x^{\varepsilon_1}< \xi\leqslant x$, and hence $\xi^{\varepsilon_1}>x^{\varepsilon}$, we may apply  \ref{cascv} to $(\xi,s_1,s_1,2\gb,\varepsilon_1,2\delta_1)$ in place of $(x,f,r,\varrho,\varepsilon,\delta)$: indeed, we check that $2\delta_1\leqslant 2\delta_0(\gb)=\ft23\beta\gb$ for $\beta$  defined by \eqref{defpar} with $\varrho=2\gb$, and $\varepsilon_1=\sqrt{\varepsilon}\geqslant 1/\sqrt{\log \xi}$. This yields
$$M(\xi;s_1)\asymp{\xi\e^{Z(\xi;s_1)}\over \log \xi}\qquad (x^{\varepsilon_1}< \xi\leqslant x).$$
Moreover, \ref{fr} immediately yields $M(\xi;r)\asymp M(\xi;r_1)$ for the same values of $\xi$.
 Therefore, for $x^{2\varepsilon_1}<z\leqslant x$, we have
$$\eqalign{M(z;r)&\gg M(z;r_1)= \sum_{m\leqslant z}t_1(m)M\Big({z\over m};s_1\Big)\cr&
\gg\sum_{m\leqslant \sqrt{z}}{zt_1(m)\e^{Z(z/m;s_1)}\over m\log (2z/m)}\asymp {z\e^{Z(z;s_1)}\over \log z}\sum_{m\leqslant \sqrt{z}}{t_1(m)\over m}\cdot\cr}$$
Now a standard manipulation resting on Rankin's method---see \eg\ \citeplus{Te16}{$(1{\cdot4})$}---yields
$$\sum_{m\leqslant \sqrt{z}}{t_1(m)\over m}\asymp \e^{Z(z;t_1)}.$$
Appealing to the identity $Z(z;s_1)+Z(z;t_1)=Z(z;r)$ completes the proof of \eqref{Mxras}.
\qed
\bigskip
We  now embark on the proof of \eqref{Mf/Mr-g}.
Without loss of generality, we may assume $\varepsilon$ arbitrarily small. Indeed, in the opposite circumstance the required estimate reduces to $M(x;f)\ll M(x;r)$ and hence follows from the inequality $|f|\leqslant r$.\par 
 For $x^{2\varepsilon_1}<z\leqslant x$, and hence $z^{\varepsilon_1}>x^\varepsilon$, we have, by the assumptions of \ref{thWg},
$$\leqalignno{
\sum_{z^{\varepsilon_1}<p\leqslant y}{\{r(p)-\re f_\tau(p)\}^\gh\log p\over p}&\ll\varepsilon_1^{2\delta_1\gh}\log y\qquad (z^{\varepsilon_1}<y\leqslant z),&(1{\cdot}4)'\cr
\sum_{z^{\varepsilon_1}<p\leqslant z}{r(p)-\re f_\tau(p)\over p}&\ll\varepsilon_1^{2\delta_1}, &(1{\cdot}7)'\cr}$$
the latter being only necessary if the extra hypothesis on $\min_{x^\varepsilon<p\leqslant x}r(p)$ is not assumed.
Therefore, we may apply \ref{fr} with $(z,r,f_\tau,\varepsilon_1,2\delta_1)$ in place of $(x,r,f,\varepsilon,\delta_1)$.
Let $L_r(x;f)$ denote the product appearing on  the right-hand side of \eqref{Mf/Mr}. Applying  this formula for $f_\tau$  while taking \eqref{Mxras} into account, we see that
$$M(z;f_\tau)=M(z;r)L_r(z;f_\tau)+O\big(\varepsilon^\delta M(z;r)\big)\qquad (x^{2\varepsilon_1}< z\leqslant x).$$
Hence
$$\eqalign{M(x;f)&=\int_1^xz^{i\tau}\d M(z;f_\tau)=x^{i\tau}M(x;f_\tau)-i\tau\int_1^xz^{i\tau-1}M(z;f_\tau)\d z\cr
&=x^{i\tau}M(x;r)L_r(x;f_\tau)-i\tau\int_{\varepsilon x}^xz^{i\tau-1}M(z;r)L_r(z;f_\tau)\d z+O\big(\varepsilon^\delta M(x;r)\big),\cr}\eqdef{Mfr}$$
where we  used, for both $\xi=\varepsilon x$ and $\xi=x$, the bound
$$\int_1^\xi{M(z;r)\over z}\d z\ll \int_1^\xi{\e^{Z(z;r)}\over \log 2z}\d z\ll{\xi\e^{Z(\xi;r)}\over \log 2\xi}\ll M(\xi;r).$$
Now, observe that, by $(1{\cdot}7)'$, we have
$$L_r(z;f_\tau)=L_r(x;f_\tau)\big\{1+O\big(\varepsilon^{\delta_1}\big)\big\}\qquad (x^{2\varepsilon_1}< z\leqslant x).$$
Since $\delta=w_f\delta_1\leqslant \delta_1$, we infer that
$$\int_{\varepsilon x}^xz^{i\tau-1}M(z;r)L_r(z;f_\tau)\d z=L_r(x;f_\tau)I(x)+O\big(\varepsilon^\delta M(x;r)\big).\eqdef{intML}$$
with
$$I(x):=\int_{\varepsilon x}^xz^{i\tau-1}M(z;r)\d z.\eqdef{Ix}$$
\par 
At this stage, we extend $s$ as defined on the primes in \eqref{defsp} to an exponentially multiplicative function by setting $s(\pnu):=s(p)^\nu/\nu!$ $(\nu\geqslant 0)$ for all primes $p$, and define the multiplicative function $t$ by the convolution formula $s*t=r$.
For $\varepsilon x\leqslant z\leqslant x$, we write
$$M(z;r)=\sum_{m\leqslant z}t(m)M\Big({z\over m};s\Big)=S+R,\eqdef{decMyr}$$
 where $S$ corresponds to the contribution of $m\leqslant \varepsilon z/x^{\varepsilon_1}$ and $R$ denotes the complementary sum.\goodbreak 
Now, observing that $M(z/m;s)=1$ for $z/x^{\varepsilon_2}<m\leqslant z$, we have
 $$\eqalign{R&\ll\sum_{\varepsilon z/x^{\varepsilon_1}<m\leqslant z/x^{\varepsilon_2}}{z|t(m)|\e^{Z(z/m;s)}\over m\log (z/m)}+\sum_{z/x^{\varepsilon_2}<m\leqslant z}{|t(m)|}\cr&\ll
 \sum_{\varepsilon z/x^{\varepsilon_1}<m\leqslant z/x^{\varepsilon_2}}{z|t(m)|\{\log (z/m)\}^{2\gb-1}\over m(\varepsilon_2\log x)^{2\gb}}+ {z\e^{Z(z;t)}\over \log z}\cr&\ll\int_{\varepsilon z/x^{\varepsilon_1}}^{z/x^{\varepsilon_2}}{z\{\log (z/v)\}^{2\gb-1}\over v(\varepsilon_2\log x)^{2\gb}}\d O\Big({v\e^{Z(z;t)}\over \log z}\Big)+{z\e^{Z(z;t)}\over \log z}\ll{z(\varepsilon_1/\varepsilon_2)^{2\gb}\e^{Z(z;t)}\over \log z}\cr&\ll {z(\varepsilon_1/\varepsilon_2)^{2\gb}\e^{Z(x;r)-Z(z;s)}\over \log z}\ll M(z;r)\varepsilon_1^{2\gb}\ll \varepsilon^\delta M(z;r),\cr}\eqdef{majR}$$
 where we used \eqref{defsp} in the form 
 $$Z(z;s)=2\gb\log (1/\varepsilon_2)+O(1)\qquad (\sqrt{x}\leqslant z\leqslant x),\eqdef{Zys}$$  
 together with the inequalities $\varepsilon_2^{2\gb}=\varepsilon^{3\gb}\leqslant \varepsilon^{9\delta_1}\leqslant \varepsilon^{\delta}$.
 \par 
 In the range $m\leqslant \varepsilon z/x^{\varepsilon_1}$, the bound~\eqref{regs} is a sufficient hypothesis to evaluate $M(z/m;s)$ by  \ref{cascv} with $(z/m,s,s,\varepsilon_1,2\delta_1)$ in place of $(x,r,f,\varepsilon,\delta)$. This furnishes the estimate
 $$M\Big({z\over m};s\Big)={\big\{1+O\big(\varepsilon^\delta\big)\big\}\e^{-2\gb \gamma}z\over \Gamma(2\gb) m(\varepsilon_2\log x)^{2\gb}}\Big(\log {z\over m}\Big)^{2\gb-1}\qquad \big(m\leqslant \varepsilon z/x^{\varepsilon_1},\,\varepsilon x\leqslant z\leqslant x\big),\eqdef{Mys}$$
 where we took into account the fact that, by \eqref{defsp}, we have $Z(x^{\varepsilon_2};s)=0$.
\par \goodbreak
Therefore, in view of \eqref{Ix}, \eqref{decMyr} and \eqref{majR}, we get
$$I(x)={\e^{-2\gb \gamma}\over \Gamma(2\gb)}J(x)+O\big(\varepsilon^\delta M(x;r)\big),$$
with
$$\eqalign{J(x)&:=\int_{\varepsilon x}^x\sum_{m\leqslant \varepsilon z^{1-\varepsilon_1}}{t(m){z^{i\tau}}\{\log (z/m)\}^{2\gb-1}\over m(\varepsilon_2\log x)^{2\gb}}\big\{1+O\big(\varepsilon^\delta\big)\big\}\d z\cr
&=\sum_{m\leqslant \varepsilon x^{1-\varepsilon_1}}\int_{\varepsilon x}^x{t(m){z^{i\tau}}\{\log (z/m)\}^{2\gb-1}\over m(\varepsilon_2\log x)^{2\gb}}\d z+R_1+R_2+R_3,\cr}\eqdef{Jx}$$
and
$$\eqalign{R_1&\ll\int_{\varepsilon x}^x \sum_{\varepsilon z^{1-\varepsilon_1}<m\leqslant \varepsilon x^{1-\varepsilon_1}}
{|t(m)|(\log z/m)^{2\gb-1}\over m(\varepsilon_2\log x)^{2\gb}} \d z\cr
&\ll {(\varepsilon_1\log x)^{2\gb-1}\over (\varepsilon_2\log x)^{2\gb}}\int_{\varepsilon x}^x\sum_{\varepsilon z^{1-\varepsilon_1}<m\leqslant \varepsilon x^{1-\varepsilon_1}}{|t(m)|\over m}\d z
\ll \sum_{\varepsilon^2 x^{1-\varepsilon_1}<m\leqslant \varepsilon x^{1-\varepsilon_1}}{x|t(m)|(\varepsilon_1/\varepsilon_2)^{2\gb}\over m\varepsilon_1\log x}
\cr
&\ll {x(\varepsilon_1/\varepsilon_2)^{2\gb}\over\varepsilon_1\log x}\int_{\varepsilon^2x^{1-\varepsilon_1}}^{\varepsilon x^{1-\varepsilon_1}}{1\over v}\d O\bigg({v\e^{Z(x;t)}\over \log x}\bigg)\ll {x(\varepsilon_1/\varepsilon_2)^{2\gb}\e^{Z(x;t)}\log (1/\varepsilon)\over\varepsilon_1(\log x)^2}\cr&
\ll {\varepsilon_1^{2\gb}\log (1/\varepsilon)x\e^{Z(x;r)}\over\varepsilon_1(\log x)^2}\ll {\varepsilon_1^{2\gb}\log (1/\varepsilon)x\e^{Z(x;r)}\over \log x}\ll
\varepsilon^\delta M(x;r),\cr
R_2&:=\sum_{\sqrt{x}<m\leqslant \varepsilon x^{1-\varepsilon_1}}{\varepsilon^\delta x|t(m)|\{\log (x/m)\}^{2\gb-1}\over m(\varepsilon_2\log x)^{2\gb}}\ll{\varepsilon^\delta x\over (\varepsilon_2\log x)^{2\gb}}\int_{\sqrt{x}}^{x^{1-\varepsilon_2}}{(\log x/v)^{2\gb-1}\over v}\d O\bigg({v\e^{Z(x;t)}\over \log x}\bigg)\cr
&\ll{\varepsilon^\delta x\e^{Z(x;t)}\over \varepsilon_2^{2\gb}\log x}\ll\varepsilon^\delta M(x;r),\cr
R_3&\ll\sum_{m\leqslant\sqrt{x}}{\varepsilon^\delta x|t(m)|\over m\varepsilon_2^{2\gb}\log x}\ll {\varepsilon^\delta x\e^{Z(x;t)}\over \varepsilon_2^{2\gb}\log x}={\varepsilon^\delta x\e^{Z(x;r)-Z(x;s)}\over \varepsilon_2^{2\gb}\log x}\asymp \varepsilon^\delta M(x;r).\cr
}$$
\par 
To evaluate the main term appearing in the right-hand side of \eqref{Jx}, it is sufficient to observe that, by partial summation, we have, for $m\leqslant \varepsilon x^{1-\varepsilon_1},$
$$\eqalign{\int_{\varepsilon x}^x{z^{i\tau}}\Big(\log {z\over m}\Big)^{2\gb-1}\d z&=\{1+O(\varepsilon)\}{x^{1+i\tau}\over 1+i\tau}\Big(\log {x\over m}\Big)^{2\gb-1}+O\Big(x\Big(\log {x\over m}\Big)^{2\gb-2}\Big)\cr&=\big\{1+O\big(\sqrt{\varepsilon}\big)\big\}{x^{1+i\tau}\over 1+i\tau}\Big(\log {x\over m}\Big)^{2\gb-1}.\cr}$$
Thus we can summarise our estimates as
$$\eqalign{I(x)&=\sum_{m\leqslant \varepsilon x^{1-\varepsilon_1}}{\e^{-2\gb\gamma}t(m)x^{1+i\tau}\{\log (x/m)\}^{2\gb-1}\over \Gamma(2\gb)(1+i\tau)m(\varepsilon_2\log x)^{2\gb}}+O\Big(\varepsilon^\delta M(x;r)\Big)\cr
&=\sum_{m\leqslant \varepsilon x^{1-\varepsilon_1}}{x^{i\tau}t(m)\over 1+i\tau}M\Big({x\over m};s\Big)+O\Big(\varepsilon^\delta M(x;r)\Big)={x^{i\tau}M(x;r)\over 1+i\tau}+O\Big(\varepsilon^\delta M(x;r)\Big)\cr}$$
Carrying this back into \eqref{intML} and \eqref{Mfr}, we obtain
$$\eqalign{M(x;f)&={x^{i\tau}L_r(x;f_\tau)M(x;r)\over 1+i\tau}+O\Big(\varepsilon^\delta M(x;r)\Big),\cr}$$
which coincides with \eqref{Mf/Mr-g}.
\medskip
\medskip
\paraun{Moments}
Given a non-negative multiplicative function $r\in\M(A,B)$ and a real additive function~$h$, let us consider the distribution function $z\mapsto F_x(z;h,r)$ of the random variable $h(n)$ on the set of integers not exceeding $x$ equipped with the measure   attributing to each integer~$n\leqslant x$ the weight $r(n)/M(x;r)$, viz.
$$F_x(z;h,r):={1\over M(x;r)}\sum_{\di{n\leqslant x}{h(n)\leqslant z}}r(n)\qquad (z\in\r).$$
\par 
Put
$$E_h(x;r):=\sum_{p\leqslant x}{r(p)h(p)\over p},\quad D_h(x;r)^2:=\sum_{p\leqslant x}{r(p)h(p)^2\over p},$$
and denote by
$\Phi(z):=\int_{-\infty}^z\e^{-u^2/2}\d u/\sqrt{2\pi}$
the distribution function of the normal law.\par 
\smallskip
The following theorem is established in \citer{Te17} as a corollary to \ref{fr}. We write $$\mu_x=\mu(x;h,r):=\max_{p\leqslant x}{|h(p)|\over D_h(x;r)},\qquad \vartheta_x=\vartheta(x;h,r):=\mu_x+1/D_h(x;r),$$
and note that $\vartheta_x\asymp \mu_x$ if $\max_{p\leqslant x}|h(p)|\gg1$ and that $\vartheta_x=o(1)$ if $D_h(x;r)\to\infty$.
\vskip-2mm\goodbreak
\propt{repad}{ (\citer{Te17})}{Let $A$, $B$,  denote positive constants. Let $x\geqslant 2$, $r\in\M(A,B)$, and let $h$ be a real additive function. Assume that:\smallskip \qquad {\rm(i)} \quad$\dsp\min_{\exp\sqrt{\log x}<p\leqslant x}r(p)\gg1$
 \ ;\qquad {\rm(ii)} \quad $D_h(x;r)\gg1$;\par \smallskip
\qquad {\rm(iii)}\quad $\mu_x\leqslant 1$;
\qquad {\rm(iv)}\quad $\dsp\sum_{\pnu\leqslant x}\sum_{\nu\geqslant 2}{r(\pnu)|h(\pnu)|\log \pnu\over \pnu}\ll1.$\note{Due to a misprint, the factor $\log \pnu$ is missing in the corresponding hypothesis of \citeplus{Te17}{cor.\thinspace2.5}.}
\par 
Then 
$$F_x\Big(E_h(x;r)+zD_h(x;r);h,r\Big)=\Phi(z)+O\big(\vartheta_x\big).\eqdef{apprep}$$
}
\par
In this section, we apply the above result to evaluate, as $x\to\infty$, the weighted moments 
$$G_m(x;r,h):={1\over M(x;r)}\sum_{n\leqslant x}r(n)\big\{h(n)-E_h(x;r)\big\}^m\qquad (m\geqslant 1).$$
We denote by $\nu_m$ the $m$th integral moment of the normal law. 
\par 
 \goodbreak
 \Propt{thmom}{Let $m\geqslant 1$ and let $h$ be a real strongly additive function. Under the hypotheses of~\ref{repad} and assuming $\vartheta_x\leqslant \dm$, we have
 $$G_m(x;r,h)=\Big\{\nu_m+O\Big(\vartheta_x(\log 1/\vartheta_x)^{m/2}\Big)\Big\}D_h(x;r)^m.\eqdef{evalGm}$$
}
\rem The assumption that $h$ is strongly additive is not essential here but it simplifies the analysis. It could be relaxed by writing $h=h_1+h_2$, where $h_1$ (resp. $h_2$) is supported on the set of squarefree (resp. squareful) integers and making {\it ad hoc} hypotheses on the values $h(\pnu)$ for $\nu\geqslant 2$. 
\medskip
\nid To lighten the writing, put $E:=E_h(x;r)$, $D:=D_h(x;r)$. With the aim of applying \eqref{apprep}, we need an upper bound for the contribution to the left-hand side of \eqref{evalGm} of large values of $|h(n)-E|$.  For $\sigma\in\r$, $|\sigma|\leqslant 1/\mu_x$, Shiu's bound \citer{Sh80} furnishes, in view of \eqref{Mxras},
$$\eqalign{\sum_{n\leqslant x}r(n)\e^{\sigma h(n)/D}&\ll x\prod_{p\leqslant x}\Big(1+{r(p)\e^{\sigma h(p)/D}-1\over p}\Big)\ll M(x;r)\exp\Big\{\sum_{p\leqslant x}{r(p)\{\e^{\sigma h(p)/D}-1\}\over p}\Big\}\cr
&\ll M(x;r)\exp\bigg\{\sum_{p\leqslant x}{\sigma r(p)h(p)\over pD}+\sum_{p\leqslant x}{\sigma^2 r(p)h(p)^2\over pD^2}\bigg\}=M(x;r)\e^{\sigma E/D+\sigma^2},\cr}$$
where we used the inequality $\e^v-1\leqslant v+v^2$ $(|v|\leqslant 1)$.
Applying this for $\sigma=\pm t,$ $0\leqslant t\leqslant 1/\mu_x$, we~get
$$\sum_{n\leqslant x}r(n)\e^{t|h(n)-E|/D}\ll \e^{t^2} M(x;r),$$
whence, selecting $t:=\sqrt{\log 1/\vartheta_x}$ and applying the above for $2t$,
$$\eqalign{\sum_{\di{n\leqslant x}{|h(n)-E|>5tD}}r(n)|h(n)-E|^m&\leqslant \sum_{n\leqslant x}r(n) {m!D^m\over t^m}\e^{2t|h(n)-E|/D-5t^2}\cr&\ll M(x;r)D^m\e^{-t^2}=\vartheta_x D^mM(x;r).\cr}\eqdef{contgdsval}$$
Now combining \eqref{contgdsval} and \eqref{apprep} yields
$$\eqalign{G_m(x;r,h)&=D^m\int_{-5t}^{5t}z^m\d\Big\{\Phi(z)+O\big(\vartheta_x\big)\Big\}+O\big(\vartheta_xD^m\big)\cr
&=\Big\{\nu_m-{1\over \sqrt{2\pi}}\int_{|z|>5t}\e^{-z^2/2}\d z+O\big(t^m\vartheta_x\big)\Big\}D^m=\Big\{\nu_m+O\big(t^m\vartheta_x\big)\Big\}D^m,\cr}$$
as required.
\qed
\bigskip
\medskip
\paraun{Sifted mean-values}
Given an arithmetic function $f$ and an integer $D$, put $f_D(n):=\1_{\{(n,D)=1\}}f(n)$. Moreover, if $f$ satisfies $\sum_{\nu\geqslant 0} |f(\pnu)|/\pnu <\infty$ for all primes $p$, define 
$$W_f(n):=\prod_{p|n}\sum_{\nu\geqslant 0}{f(p^\nu)\over p^\nu}\qquad (n\geqslant 1),$$
and denote by $P^+(n)$ the largest prime factor of an integer $n>1$, with the convention that $P^+(1)=1$.
\par 
In his paper \citer{El17}, Elliott indicates that the following effective estimate is Çprovided by combining the argument of Elliott \& Kish \citer{EK16b},  \citeplus{El17}{th. 2}, with that of  the  taxonomy section of Elliott \& Kish~\citer{EK16a}È:
{\it Let $r$ denote a non-negative, exponentially  multiplicative function, satisfying, for suitable positive constants $A$, $\gb$, and $c_1$,
$$\sup_pr(p)\leqslant A, \quad\sum_{w<p\leqslant v}{r(p)-\gb\over p}\geqslant -c_1\qquad (\ft32\leqslant w\leqslant v).\eqdef{cnE}$$
Then, uniformly for  $x\geqslant 2$, $D\geqslant 1$, $P^+(D)\leqslant x$, we have
$$M(x;r_D)=M(x;r)\bigg\{{1\over W_r(D)}+O\bigg({(\log_22D)^{1+A}\over (\log x)^{\gc}}\bigg)\Bigg\}\eqdef{faPE}$$
where $\gc:=\gb^3/\{\gb^2+3456 A^2\}$ provided $\gb\leqslant 12\sqrt{2A}$.\note{The condition $P^+(D)\leqslant x$ is omitted in \citer{El17}. We reinserted it  since it does not involve any loss of generality.}}
\par \goodbreak\smallskip
The above statement can be directly compared with an almost immediate consequence of \ref{fr}. Put
$$\beta=\beta(\gb,A):=1-{\sin(\pi\gb/A)\over \pi\gb/A},\qquad  \delta=\delta(\gb,A):=\ft1{12}\gb\beta,\eqdef{delta0+}$$
and recall definition \eqref{C0} of the class $\M(A,B)$. Note that, if  $\gb\leqslant A$, we have $\beta\geqslant \gb^2/A^2$ by the product formula for $\sin(\pi z)/\pi z$, and hence $\delta\geqslant \gb^3/12 A^2$ .\par 
\Propt{gD_GT}{Let $A>0$, $B>0$, $r\in \M(A,B)$, $r\geqslant 0$, $x\geqslant e^4$. Assume that, for a suitable constant~$\gb$, $0<\gb\leqslant \min(1,A)$, and  $\eta_x:=(\log x)^{-1/4}$, we have
$$\sum_{y<p\leqslant y^{1+\eta_x}}{r(p)\log p\over p}\geqslant \gb\eta_x\log y\qquad \big(\e^{1/\eta_x}\leqslant y\leqslant x^{1/(1+\eta_x)}\big).\eqdef{moygp}$$
Then,  uniformly for $x\geqslant 2$, $D\geqslant 1$, $P^+(D)\leqslant x$, we have
$$M(x;r_D)=M(x;r)\bigg\{{1+O(\chi)\over W_r(D)}+O\Big({1\over (\log x)^{\delta/2}}\Big)\bigg\},\eqdef{faGT}$$ where
$\chi:=\1_{\{\log_23D>(\log x)^{\gb^3/17A^3}\}}$.}
\par \medskip
\rems (i) The restriction to exponentially multiplicative functions has been dropped.\par 
(ii) Condition \eqref{moygp} is significantly less restrictive than \eqref{cnE}.
\par 
(iii) The error term of \eqref{faGT} is always smaller than that of \eqref{faPE} by a power of $\log x$.
\medskip
\nid We shall  apply \ref{fr} to the pair $(r,f)=(r,r_D)$, replacing $(A,\gb)$ by $(A/2,\gb/4)$, and selecting $\varepsilon:=1/\sqrt{\log x}$, so that $\varepsilon_1=\eta_x$. \par 
First consider the case when $\log_23D\leqslant (\log x)^{\gb^3/17A^3}$.  We then have, for large $x$,
$$\eqalign{\sum_{p|D}{r(p)\over p}\leqslant A\log_33D+O(1)\leqslant {\gb^3\over 16A^2}\log_2x\leqslant \ft18\beta\gb\log (1/\varepsilon).\cr}$$
Therefore, condition \eqref{pqr} holds for our modified parameters. To check that conditions \eqref{cdpr1} with $\gh:=(4-\gb)/\gb$, and  \eqref{cdpr} with $\gh=1$ are also fulfilled, we note that a well known estimate provides
 $$\sum_{p|D}{\log p\over p}\ll \log_23D.$$
We hence have
$$\sum_{\di{x^\varepsilon<p\leqslant y}{p|D}}{\log p\over p}\ll (\log x)^{\gb^3/16A^3}\leqslant {\log y\over 
(\log x)^{\frac12(1-\gb^3/8A^3)}}\leqslant {\log y\over (\log x)^{7/16}}\qquad (x^\varepsilon<y\leqslant x).$$
Since, for $\gh:=(4-\gb)/\gb$, we have $\ft7{16}>\dm\delta\gh=\ft1{24}\beta(4-\gb)$, we see that \eqref{cdpr1}   holds for this value of~$\gh$ and $\delta_1=\delta$.
\par 
Next,
$$\sum_{\di{x^\varepsilon<p\leqslant x}{p|D}}{1\over p}\leqslant \sum_{p|D}{\log p\over p\sqrt{\log x}}\ll {1\over (\log x)^{\frac12(1-\gb^3/8A^3)}}\ll\varepsilon^{\delta},$$
with a lot to spare, since $\delta\leqslant \frac1{12}$, $1-\gb^3/8A^3\geqslant \frac78$.  This shows that \eqref{cdpr} with $\gh=1$ is satisfied.
\par\smallskip Applying \ref{fr} with the parameters defined above,  we get in the case under consideration
$$M(x;r_D)=M(x;r)\Big\{{1\over W_r(D)}+O\Big({1\over (\log x)^{\delta/2}}\Big)\Big\},$$
which is compatible with \eqref{faGT}.
\par 
If $\log_23D> (\log x)^{\gb^3/17A^3}$,  estimate \eqref{faGT} reduces to the Halberstam-Richert upper bound \citer{HR79} since, by \ref{Mzr}, $M(x;r)\asymp x\e^{Z(x;r)}/\log x$.
\qed
\medskip
{\noi\bf Acknowledgement.} The author thanks RŽgis de la Bretche for useful comments on a preliminary draft of this article. 
\goodbreak
\par\medskip
\bigskip\bigskip
\centerline{\twelvebf References}\bigskip
{\leftskip5mm\rightskip5mm\eightpoint{
\bibtem{El17} P.D.T.A Elliott, Multiplicative function mean values: asymptotic estimates, {\it Funct. Approx. Comment. Math. \bf56},\numero2 (2017),  217Ð238. \par 
\bibtem{EK16a}  P.D.T.A Elliott \& J. Kish, Harmonic analysis on the positive rationals I: Basic results, {\it J. Math. Sci. Univ. Tokyo \bf23},\numero3 (2016),  569Ð614.\par 
\bibtem{EK16b} P.D.T.A Elliott \& J. Kish, Harmonic analysis on the positive rationals II: Multiplicative functions and Maass forms, {\it J. Math. Sci. Univ. Tokyo \bf23},\numero3 (2016), 615Ð658. \par 
\bibtem{Ha68} G. Hal‡sz, †ber die Mittelwerte multiplikativer zahlentheoretischer Funktionen,
{\it Acad. Math. Acad. Sci. Hungar. \bf19} (1968), 365--403.\par 
\bibtem{HR79} H. Halberstam \& H.-E. Richert, On a result of R.R. Hall, {\it J. Number Theory \rm (1) \bf 11}
(1979), 76--89.\par 
\bibtem{IK25} K.-H. Indlekofer \& E. Kaya, Estimates for multiplicative functions, II, {\it Annales Univ. Sci. Budapest., Sect. Comp. \bf 58} (2025), 1-11.
\bibtem{Sh80}
P. Shiu, A Brun-Titchmarsh theorem for multiplicative functions, {\it J. Reine Angew. Math.} {\bf 313} (1980), 161--170.
\par 
\bibtem{Te16} G. Tenenbaum, Fonctions multiplicatives, sommes d'exponentielles, et loi des grands nombres, {\it Indag. Math. \bf27} (2016), 590--600.
\bibtem{Te17} G. Tenenbaum, Valeurs moyennes effectives de fonctions multiplicatives complexes, {\it Ramanujan J. \bf44}, \numero3 (2017), 641-701; Corrig. {\it ibid. \bf 51}, \numero1 (2020), 243-244.\par
\bibtem{Wi67} E. Wirsing, Das asymptotische Verhalten von Summen \"uber multiplikative Funktionen
II, {\it Acta Math. Acad. Sci. Hung. \bf 18} (1967), 411--467.

\par }
}
\vskip 5mm
{\sevenrm\baselineskip9pt
G\'erald Tenenbaum\par
Institut \'Elie Cartan\par 
Universit\'e de Lorraine\par
 BP 70239\par
54506 Vand\oe uvre-ls-Nancy Cedex\par
 France
\smallskip
e-mail : \seventt gerald.tenenbaum@univ-lorraine.fr\par}
\end